\documentclass[11pt]{article}
\usepackage{amssymb}
\usepackage{amsmath, amscd}

\newtheorem{theorem}{Theorem}[section]
\newtheorem{rem}[theorem]{Remark}
\newtheorem{defi}[theorem]{Definition}
\newtheorem{lemma}[theorem]{Lemma}

\newtheorem{prop}[theorem]{Proposition}
\newtheorem{conj}[theorem]{Conjecture}

\def \C{{\mathbb C}}
\def \R{{\mathbb R}}
\def \Z{{\mathbb Z}}
\def \P{{\mathbb P}}
\def \Q{{\mathbb Q}}
\def \H{{\mathbb H}}

\newcommand{\Teich}{\operatorname{\sf Teich}}

\newcommand{\Comp}{\operatorname{\sf Comp}}

\newcommand{\Hyp}{\operatorname{\sf Hyp}}
\newcommand{\Diff}{\operatorname{Diff}}
\newcommand{\Per}{\operatorname{\sf Per}}
\newcommand{\Pos}{\operatorname{\sf Pos}}
\newcommand{\Nef}{\operatorname{\sf Nef}}
\newcommand{\Def}{\operatorname{\sf Def}}
\newcommand{\Bir}{\operatorname{\sf Bir}}
\newcommand{\Kah}{\operatorname{\sf Kah}}
\newcommand{\Stab}{\operatorname{Stab}}
\newcommand{\Mov}{\operatorname{\sf Mov}}
\newcommand{\Mon}{\operatorname{Mon}}
\newcommand{\Hdg}{\operatorname{\sf Hdg}}
\newcommand{\Gr}{\operatorname{\sf Gr}}

\newcommand{\Aut}{\operatorname{Aut}}
\newcommand{\Tw}{\operatorname{Tw}}

\begin{document}

\title{Rational curves and MBM classes on hyperk\"ahler manifolds: a survey} 
\author{E. Amerik, M. Verbitsky}
\date{\today}
\maketitle

\section{Introduction}

This paper deals with rational curves and birational
contractions on hyper\-k\"ahler manifolds. Here
``hyperk\"ahler'' means compact K\"ahler holomorphically 
symp\-lec\-tic manifold. We survey
some recent results about minimal rational curves,
their deformations, extremal rays associated with these
curves, and the geometry of the K\"ahler cone.

It is known since Shafarevich and Pyatetski-Shapiro's work
\cite{PS-S} in the algebraic case and Looijenga-Peters
\cite{LP} in the compact K\"ahler case that some important
aspects of the geometry of K3 surfaces are governed by
smooth rational curves on these surfaces. By adjunction,
the square of such a curve is always equal to $-2$. One
calls it a $(-2)$-curve. A fundamental theorem directly
concerned with the subject of this paper states that a
line bundle on a K3 surface is ample if and only if it is
of positive square and positive on all $(-2)$-curves. We
also have an analogue in the non-algebraic case, with an
ample line bundle replaced by a K\"ahler class. In other
words, the orthogonal hyperplanes to the $(-2)$-curves bound
the ample (resp. K\"ahler) cone inside the positive cone
in $NS(X)\otimes {\R}$ (resp. $H^{1,1}_{\R}(X)$).

The number of $(-2)$-curves on a K3 surface $S$ is not necessarily
finite. However H. Sterk \cite{St} proved that there
are finitely many of them up to the action of $\Aut(S)$. In the
projective case he actually proved a more precise result:
the group $\Aut(S)$ has a polyhedral fundamental domain on
the ``rational hull of the ample cone'' $\Nef^+(X)$ (which
by definition is the ample cone to which one attaches the
part of the boundary defined over the rationals).

In this survey we consider the {\bf irreducible holomorphic symplectic manifolds (IHSM)}, also known as {\bf simple hyperk\"ahler manifolds}, which are higher-dimensional analogues of K3 surfaces. 

\begin{defi} An irreducible holomorphic symplectic manifold is a simply-connected compact K\"ahler manifold $X$ whose space of holomorphic two-forms $H^{2,0}(X)$ is generated by a symplectic (that is, nowhere degenerate) form $\sigma$. 
\end{defi}
	
In particular such an $X$ has even dimension $2n$ and trivial canonical class. There are, at present, two infinite series of families and two ``sporadic'' families of IHSM known, and so far we are unable to answer the question whether there are other examples. 

\medskip

{\bf Example:} If $S$ is a K3 surface, then the $n$-th punctual Hilbert scheme $S^{[n]}$ is an IHSM, see \cite{Be}. The deformations of $S^{[n]}$ are, at this date, the most studied hyperk\"ahler manifolds. One calls them IHSM {\bf of K3 type}.

\medskip

 The analogy with a surface is especially striking because
 on the second cohomology $H^2(X, \Z)$ of such a manifold
 $X$ there is an integral quadratic form $q$, the {\bf
   Beauville-Bogomolov form} of signature
 $(1, b_2-3)$ on $H^{1,1}(X,\R)$ (and $(3, b_2-3)$ on the whole second cohomology). This
 form is defined by integrating differential forms against
 appropriate powers of $\sigma$ and $\bar \sigma$ in
 \cite{Be}, but in fact is of topological nature by
 \cite{F}.

\medskip

{\bf Example:} Let $S^{(n)}$ be the $n$-th symmetric power of a K3 surface $S$ and $HC: S^{[n]}\to S^{(n)}$ the Hilbert-Chow map.
The cohomology class of the exceptional divisor of $HC$ is divisible by 2, and the second cohomology group of $S^{[n]}$ is naturally isomorphic to the direct sum of $H^2(S,\Z)$ and $\Z e$, where $e$ is one half of the exceptional divisor. This direct sum turns out to be orthogonal 
with respect to $q$. Moreover $q$ restricts to the intersection form on $H^2(S,\Z)$, and $q(e)=-2(n-1)$. The second cohomology lattice remains the same under deformations, so this description can be used to study 
the IHSM of K3 type. The difference between
$S^{[n]}$ and its arbitrary deformation is, roughly
speaking, that the class $e$ does not have to be a Hodge
class on the latter.

\medskip

All K\"ahler classes on an IHSM $X$ have positive Beauville-Bogomolov square. The set of classes with positive square in 
$H^{1,1}(X,\R)$ has two connected components. We call {\bf the positive cone} the component containing the K\"ahler classes.
The following analogue of Shafarevich and Pyatetski-Shapiro's/Looijenga and Peters' result
is available.

\begin{theorem}\label{HB} (Huybrechts, Boucksom; 
\cite{Bou,_Huybrechts:cone_})  The K\"ahler classes on an
IHSM $X$ are the elements of the positive cone $\Pos(X)$
which are strictly positive on all rational curves on $X$.
\end{theorem}

Thus the K\"ahler cone is cut out within the positive cone by the orthogonal hyperplanes to the classes of the rational curves.
In the same way in the projective case, the ample classes are elements of the Neron-Severi group which have positive square and are positive on all the rational curves.

If we pursue the analogy with the surface case, the next natural question is: what can be said about the possible Beauville-Bogomolov squares of these curves? Note that using the Beauville-Bogomolov form, the curve classes are naturally viewed as classes in $H^2(X, \Q)$, that is why we can talk of the Beauville-Bogomolov square of such a class. It is a rational number. One might also wonder about a possible connection between the precise value of this square and the geometry of the corresponding rational curves on $X$.

Hassett and Tschinkel initiated the study of this question in \cite{HT1}, dealing with projective IHS fourfolds, especially those which are of K3 type. In the projective case, the classes we want to understand are exactly those generating the ``extremal rays'' of the Mori cone in birational geometry. Each curve class has two numerical invariants: its Beauville-Bogomolov square and its ``denominator'', that is, the image in the discriminant group $L^*/L$, where $L=H^2(X, \Z)$, and $H_2(X, \Z)$ is viewed as its dual embedded in 
$H^2(X, \Q)$, just as we have mentioned above. For IHS manifolds of K3 type of dimension $2n$, the discriminant group is cyclic of order $2n-2$. Hassett and Tschinkel have observed that in the fourfold examples one sees three types of extremal rays: integral with square $(-2)$, half-integral with square $-1/2$ and half-integral of square $-5/2$. They conjectured that these are indeed all the extremal rays of the cone of curves. Equivalently, the ample cone
consists of such classes in $NS(X)\otimes \R$ which have positive square and are positive on all such curves.

By Kawamata-Shokurov base-point-free theorem (see e.g. \cite{KM}), extremal rays on a projective IHS manifold can be contracted. Hassett and Tschinkel, describing 
several examples of such contractions, have also conjectured that one could read the geometry of the contraction locus from 
the numerical invariants of the ray. Namely a half-integer class of square  $-5/2$ should give rise to a lagrangian plane (contracted to a point) and every other class to a family of rational curves over a K3 surface (contracted to the K3 surface).

Part of these conjectures has been proved in the subsequent papers by Hassett and Tschinkel; also an effort has been made to state and prove similar conjectures in higher dimension. The initial conjecture of Hassett and Tschinkel on the numerical description of the ample cone on an IHS manifold of K3 type, however, turned out to be false. A correct formula in terms of the Mukai lattice has been given by Bayer and  Macri \cite{BM} for punctual Hilbert schemes and more generally, moduli spaces of stable sheaves on a K3 surface. The extension from a punctual Hilbert scheme to its deformation can be found in \cite{BHT}. Bayer and Macri's work 
is highly involved and is certainly out of the scope of the present survey; we shall see, 
however, that in low dimensions one can recover their results by fairly elementary means.

Concerning the finiteness property, the extension of Sterk's theorem to the higher-dimensional case is known as the Morrison-Kawamata cone conjecture and attracts considerable attention since early 1990's.

\begin{conj}\label{cone-conj} 
(\cite{_Kawamata:cone_,_Morrison:Beyond_}
Let $X$ be a projective
  manifold with trivial canonical class. Then $\Aut(X)$
  acts with finitely many orbits on the set of faces of
  its ample cone. Moreover denote by $\Nef^+(X)$ the
  ``rational hull'' of the ample cone, that is the
  smallest convex 
cone containing the ample cone and all rational points of
its boundary. Then $\Aut(X)$ has a rational polyhedral
fundamental domain on $\Nef^+(X)$.
\end{conj}

In this generality, the conjecture is surprisingly
difficult and still wide open even for Calabi-Yau
threefolds (though there are some partial results, notably
by Kawamata \cite{_Kawamata:cone_}  and Totaro \cite{_Totaro:MK_cone_}). 
Apart from Sterk's original result, little
has been known until recently. For holomorphic symplectic
manifolds, the ``birational'' version of this conjecture
has been established by Markman in \cite{Mar-survey}; later Markman and
Yoshioka in \cite{MY} have proved the cone conjecture for IHS manifolds
of K3 or generalized Kummer type. The cone conjecture for
IHS manifolds in general has been proved in
\cite{AV-cone}, \cite{AV-orbits}.

In this survey we summarize this one and some of our other contributions to the field in recent years. Our techniques are largely based on deformations, in particular to the non-algebraic manifolds which often have fairly simple geometry. Deformation theory in the hyperk\"ahler context is well-understood and many computations are easy and elementary. The definition of MBM classes arises naturally from the urge to put the notion of an extremal rational curve in the deformation-invariant context.

The crucial ingredient which made it possible to prove the
cone conjecture for IHS manifolds was the ergodic theory,
more particularly Ratner's theorems about orbits of actions
of unipotent groups on homogeneous spaces. These techniques
have been initially used by the second-named author in
order to prove the density of orbits of the mapping class group
action on the Teichm\"uller space (\cite{V-erg}), and apply successfully
to the cone conjecture setting in \cite{AV-cone}.

Unless otherwise specified, by a rational curve we mean the image of a generically injective map from $\P^1$, so that it can be 
singular, but not reducible.

\section{MBM classes: equivalent definitions and basic properties}

\subsection{Deforming rational curves: first remarks}

The Riemann-Roch formula implies that any divisor of square $-2$ on a K3 surface is effective, up to a $\pm$ sign. However not all Hodge classes of square $-2$ on a K3 surface, up to a $\pm$ sign, are represented by $(-2)$-curves. Indeed the corresponding effective divisor can be reducible. 

Nevertheless, for any $(-2)$-class $z$ on a K3 surface $X$
there is a deformation $X'$ where $H^{1,1}(X')\cap
H^2(X',\Q)=\langle z\rangle$,  or, equivalently,
the Picard group is cyclic and generated by $z$.
Indeed, the locus where
a class $\alpha$ remains of type $(1,1)$ is equal to
$\alpha^{\bot}$ in the
space of local deformations $\Def(X)$, identified with a
neighbourhood of zero in $H^1(X, \Omega^1_X)\cong
H^1(X,T_X)$, and a general element of the hyperplane $\alpha^{\bot}$ is not orthogonal to any other integral class. Also the same is true
globally in the period domain, interpreted as an open subset of
a quadric in $\P H^2(X,\C)$ (Definition \ref{_Period_Definition_}). 
The orthogonality here is taken
with respect to the intersection form; see Subsection
\ref{_Param_Subsection_}
for the relevant definitions.  
On such a deformation, $z$ (up to 
a $\pm$ sign) is necessarily represented by a smooth
rational curve, and there are no other curves on $X$.

The notion of an MBM class is inspired by this
observation, which works in the same way in higher
dimension.

Indeed, similar arguments can be applied to a 
 higher-dimensional IHSM $X$. The space
of deformations of $X$ where a class $\alpha\in H^2(X, \Z)$ is of
type $(1,1)$ can be described as the hyperplane
$\alpha^{\bot}$ (where the orthogonality is taken with respect
to the Beauville-Bogomolov form).
A general deformation $X'$ of this type satisfies
$H^{1,1}(X')\cap H^2(X',\Q)= \langle \alpha\rangle$.
We can now distinguish two cases: it can
happen that on such a generic deformation $X'$, 
some multiple of $\alpha$
is represented by a curve (recall that we view curve classes as rational $(1,1)$-classes, via the Beauville-Bogomolov form), or that no multiple of $\alpha$
is represented by a curve\footnote{As we are dealing with manifolds which are not necessarily projective, this second case is indeed possible.}. By Huybrechts-Boucksom's
theorem (Theorem \ref{HB}), in the second case every
class in the positive cone of $X'$ is K\"ahler. In the first case, a
multiple of $\alpha$ is represented by a rational curve
and $\alpha^{\bot}$ defines the unique wall of the
K\"ahler cone of $X'$. 

Let us suppose from now on that we are in the first case.
The crucial observation is that by deformation theory of rational curves in IHS manifolds, this is also the case for {\it any} deformation $X'$ with Picard group generated by $\alpha$ over the rationals. Let us briefly summarize the basic facts on the deformations of rational curves (see e.g. \cite{AV-MBM} for details).

\begin{prop}\label{def-rc} (the first part is due to Ran \cite{Ran}, Cor. 5.1, with the proof below communicated to us by Markman) Let $C\in X$ be
  a rational curve on an IHSM
  $X$ of dimension $2n$. Then 
	\begin{itemize}
		\item $C$ deforms within $X$ in a family of dimension at least $2n-2$. 
		\item If, moreover, $C$ is minimal (that is, $C$ cannot be bent-and-broken\footnote{This holds if $C$ is of minimal degree, with respect to a K\"ahler class, among the rational curves in the uniruled subvariety of $X$ covered by the deformations of $C$.}), then $C$ deforms within $X$ in a family of dimension exactly $2n-2$. 
			
		\item Finally, if ${\cal X}\to \Def(X)$ denotes the local universal family, then a minimal rational curve $C$ deforms exactly to those fibers $X_t$ that keep the cohomology class of $C$ of type $(1,1)$.
	\end{itemize}
\end{prop}

The idea of proof of the proposition is that the deformation spaces of rational curves on IHS manifolds have greater (by one) dimension than what we can expect from the Riemann-Roch formula. An extra parameter comes from the hyperk\"ahler structure in the differential-geometric setting, the {\bf twistor family}. In fact it follows from the Calabi-Yau theorem that for any 
IHS $X=(M,I)$ (where $M$ denotes the underlying differentiable manifold and $I$ the complex structure defining $X$)
there is a hyperk\"ahler metric $g$ on $M$ inducing a whole 2-sphere (or complex projective line) of complex structures
$aI+bJ+cK$, where $I,J,K$ multiply like quaternions and $a^2+b^2+c^2=1$. Every IHS $X$ thus comes as a member of a 
''twistor family'' $\Tw(X)\to \P^1$. We can consider the deformations of $C$ in $\Tw(X)$: they are the same as those of $C$ in $X$ since 
the class of $C$ is no longer Hodge in the neighbouring complex structures (most of which carry no curves at all). Since the dimension of $\Tw(X)$ is greater by one, 
we also get an extra parameter by Riemann-Roch (which involves the rank of the pull-back of the tangent bundle of the ambient manifold).

In the proof of the second assertion, the rational
quotient fibration, also known as MRC fibration, of the uniruled subvariety $Z\subset X$
covered by deformations of $C$ is involved as a main tool. 
One obtains the following proposition as a byproduct.

\begin{prop}\label{rat-quot} The codimension of $Z$ in $X$ is equal to the relative dimension of its MRC fibration and to the corank of the restriction of the symplectic form $\sigma$ to $Z$ at a general point. In particular $Z$ is coisotropic, and if $Z$ is a divisor then it is birationally a $\P^1$-bundle.
\end{prop}

To get the statement on the walls of the K\"ahler cone from proposition \ref{def-rc}, one applies it to the minimal rational curves which always 
have class proportional to $\alpha$ on a generic deformation $X'$ preserving $\alpha$ as a Hodge class: indeed on such $X'$ there are simply no other Hodge classes. One gets that the minimal rational curves deform together with their cohomology class $\alpha$ and thus define a wall of the K\"ahler cone on $X'$.

On $X$ itself, the group of Hodge classes is greater and rational curves in the class $\alpha$ are not necessarily minimal. But though $\alpha^{\bot}$ does not necessarily define a wall of the K\"ahler cone of $X$,
it turns out that it is closely related to such walls. We call $\alpha$ an MBM (``monodromy birationally minimal'') class. We shall give the precise definition  2.10, and we shall also formulate a theorem explaining this name. In order to do this we need some discussion concerning the parameter spaces for IHS manifolds. 

Before dealing with the parameter spaces in the next subsection, we recall one important source of rational curves on IHS manifolds.

\begin{theorem} (Boucksom, \cite{Bou}, Proposition 4.7) Every prime divisor with negative Beauville-Bogomolov square on an IHS manifold is uniruled.
\end{theorem}

The converse is not true in general. However in this survey we are particularly interested in the description of the Kahler cone inside the positive cone; if we want to describe its walls as orthogonals to some classes, then by Hodge index theorem these classes necessarily have negative 
square.

In Markman's terminology \cite{Mar-survey}, the prime divisors of negative Beauville-Bogomolov square are called {\bf prime exceptional divisors}.
Indeed by a result of Druel \cite{D} they can be contracted, possibly after a sequence of flops. We shall return to contractibility issues later.

\subsection{Parameter spaces for hyperk\"ahler manifolds}
\label{_Param_Subsection_}

Let $M$ be the underlying differentiable manifold of an IHS $X$: we view $X$ as a pair $(M,I)$ where $I$ is a complex structure on $M$.
We consider the space $\Comp$ of all complex structures of
K\"ahler type on $M$. By  Kodaira-Spencer stability theorem
(\cite{_Kod-Spen-AnnMath-1960_}), this
is an open subset in the space of all complex structures. The group $\Diff$ of diffeomorphisms of $M$ and its subgroup of isotopies $\Diff_0$ 
act on $\Comp$.

\begin{defi} The Teichm\"uller space $\Teich=\Teich(M)$ is $\Comp/\Diff_0$, the quotient of  $\Comp$ by isotopies.
	The mapping class group is $\Gamma=\Diff/\Diff_0$.
\end{defi}

It turns out that the Teichm\"uller space is a possibly non-Hausdorff finite-dimensional smooth complex space. 
One would like to consider the quotient $\Teich/\Gamma$ as the ``moduli space'' of IHS. Unfortunately this makes little sense. Indeed Verbitsky proved in \cite{V-erg} that most of the orbits of the action of $\Gamma$ on $\Teich$ are dense.

The space $\Teich$ also may have infinitely many connected components. This is the case exactly when the subgroup
$K\subset \Gamma$ acting trivially on $H^2(M)$, called the Torelli subgroup, is infinite. In fact this group acts on 
$\Teich$ by permuting its connected components (\cite{V-Tor-err}) and the space $\Teich/K$ already has finitely many connected components which can be identified to those of 
$\Teich$. In a more algebraic setting, $\Teich/K$ is also known as the {\bf moduli space of marked IHS manifolds}.

We denote by $\Teich^I$ the connected component of $\Teich$ containing our given complex structure $I$. 
The {\bf monodromy group} $\Mon_I$ is the subgroup of $\Gamma$ preserving $\Teich^I$.
By \cite[Corollary 7.3]{V-Tor}, its intersection with $K$ is finite. Sometimes we drop the index $I$, if no risk of confusion arises.
The {\bf Hodge monodromy group} $\Mon^{\Hdg}$ is the part of the monodromy which preserves the Hodge type.

The crucial facts about the monodromy group are as follows.

\begin{theorem} (\cite{V-Tor-err}, Theorem 2.6) $\Mon_I$ has a finite index image (as well as a finite kernel) under the natural map $\Phi$ to $O(H^2(M,\Z), q)$.
\end{theorem}

Consequently, the Hodge monodromy maps to a subgroup of finite index in $O(NS(X), q)$.

\begin{theorem} The monodromy group acts ``ergodically'' on $\Teich^I$, that is, most of its orbits are dense (namely the union of non-dense orbits is of measure zero).
\end{theorem}

More precisely, using Ratner theory on unipotent group actions on homogeneous spaces, Verbitsky gave a classification of orbit closures of this action. In order to explain it we first need to introduce the period map for hyperk\"ahler manifolds.

\begin{defi} \label{_Period_Definition_}
The {\bf period domain} is the space 
\[ {\mathcal P}=\{l \in \P H^2(M,\C)| q(l)=0,\ q(l,\bar l)>0\}. 
\]
The {\bf period map} $\Per: \Teich\to {\mathcal P}$
assigns $H^{2,0}(X)=H^{2,0}(M,I)$ to a complex structure
$I$.
\end{defi}

Verbitsky proved in \cite{V-Tor} that the period map is an isomorphism on each connected component $\Teich^I$ of 
$\Teich$,
provided that one glues together the non-separated points of $\Teich^I$ ({\bf a global Torelli theorem for hyperk\"ahler manifolds}). Some more precise versions of this result 
have later been given by Markman in \cite{Mar-survey}. In particular Markman gives the following description of the non-separated points glued together by the period map.

For a complex IHS manifold $X=(M,I)$ the space $H^{1,1}(X,\R)$ only depends on the period point. The positive cone
$\Pos(X)\subset H^{1,1}(X,\R)$ has a chamber decomposition in the {\bf birational K\"ahler cone}\footnote{Our terminology here is non-standard, to save space and keep notations simple. Markman calls this chamber the {\bf fundamental exceptional chamber} whereas by the birational K\"ahler cone of $X$ one usually means the union of the K\"ahler cones of the IHS birational models of $X$, i.e. the complement of a certain number of walls in the fundamental exceptional chamber. We prefer not to separate the two notions and add the ``interior'' walls to the union of the K\"ahler cones of the birational models.} $\Bir(X)$ and its monodromy transforms $\gamma(\Bir(X))$, $\gamma\in \Mon$. The walls of this decomposition are given precisely by hyperplanes orthogonal to the classes of prime divisors of negative Beauville-Bogomolov square on $X$ (prime exceptional divisors; see end of Section 2.1). 

The chamber $\Bir(X)$ is in its turn subdivided into K\"ahler cones of the birational models of $X$, and analogously for $\gamma(\Bir(X))$. The resulting smaller chambers are called the K\"ahler chambers.

\begin{theorem} (Markman \cite{Mar-survey}) The points of $\Per^{-1}(\Per(I))$ are in one-to-one correspondence with the K\"ahler chambers in $\Pos(X)\subset H^{1,1}(X,\R)$.
\end{theorem}

Just in the same way as for the local deformations, the
complex structures where a given class $z\in H^2(M, \Z)$
is of type $(1,1)$ are parameterized by
$\Per^{-1}(z^{\bot})$, and their period points by the
hyperplane section $z^{\bot}\subset \mathcal P$. For a very
general complex structure $X=(M, I)$, there are no
integral (or rational) classes of type $(1,1)$. In this case the
K\"ahler cone $\Kah(X)$ is equal to the positive cone and
the Teichm\"uller space is separated at $I$. The same
holds if $I$ does have rational classes of type $(1,1)$
but they are all of non-negative Beauville-Bogomolov square,
or are not represented by curves. For other IHS the
chamber decomposition can be rather complicated, in
particular $\Teich$ is not separated even for K3
surfaces. For K3 surfaces, the chamber decomposition is
given by the orthogonals to $(-2)$-classes and $\Bir$
coinsides with $\Kah$ (see the beginning of this text for
some finiteness results in this case).

It is often useful to view $\mathcal P$ as a real variety of Grassmannian type rather than an open subset of a complex quadric. Namely by taking the real and the imaginary part of a complex line $l$ one observes that $\mathcal P$ is identified to
the space of oriented positive 2-planes in $H^2(M,\R)$:
$${\mathcal P}=Gr_{++}(H^2(M,\R))\cong SO(3, b_2-3)/SO(2)\times SO(1, b_2-3).$$ 
The {\bf monodromy orbit classification} by Verbitsky, provided
that $b_2\geq 5$,  is now as follows (\cite{V-erg-err}):

\begin{itemize} 
	\item If the period plane $P$ is rational, the 
		monodromy orbit of the complex structure is closed;

	\item If the period plane $P$ contains no rational vectors, the orbit is dense;

	\item If the period plane contains a single
          rational vector $v$ up to proportionality, the
          image of the orbit closure by the period map is
		the set of all period planes containing a monodromy image of $v$. So it has, in general, countably many irreducible components. Each component (consisting of the planes containing a fixed rational vector) is 
		a totally real submanifold of the period domain, that is, a fixed point set of an antiholomorphic involution. 
		In particular, looking at its tangent space we see that
	the orbit closure is not contained in a proper
          complex submanifold, even locally, nor contains
          any positive-dimensional complex subvarieties.

\end{itemize}

This classification is an essential ingredient for the results on contractibility and deformation invariance in Section 4.

\subsection{MBM classes}\label{mbm-def}

Let $z\in H^2(M,\Z)$ be a negative class, that is, $q(z)<0$. Denote by $\Teich_z$ the part of the Teichm\"uller space formed by the complex structures in which $z$ is a Hodge class (so that $\Teich_z=\Per^{-1}{\mathcal P}_z$, where ${\mathcal P}_z=z^{\bot}$). The following theorem is proved  in \cite{AV-MBM}.

\begin{theorem}
The following properties of $z$ are equivalent:

\begin{itemize} 
	\item $\Teich_z$ contains no twistor lines;

         \item For all $I\in \Teich_z$ such that $X=(M,I)$ has Picard rank one, $X$ contains a rational curve.

\item For generic $I\in \Teich_z$, $X=(M,I)$ contains a rational curve.

\item For $I\in \Teich_z$, $z^{\bot}$ contains a wall of a K\"ahler chamber, that is, for some $\gamma\in \Mon^{\Hdg}$, 
$\gamma(z)^{\bot}$ contains a wall of the K\"ahler cone of an IHS birational model $X'$ of $X=(M,I)$.

\end{itemize}
\end{theorem}

\begin{defi} An MBM class is a negative class in $H^2(M,\Z)$ which has these equivalent properties.
\end{defi}

If $z$ is MBM, then up to a rational multiple it is represented by a rational curve on a general $X=(M,I)$,
$I\in \Teich_z$. On a special $X$, such a curve can become reducible and the classes of the components might be not
proportional to $z$. In particular this means that $z^{\bot}$ no longer contains a wall of the K\"ahler cone, though by definition it remains MBM. The classes of the irreducible components are also MBM, provided they have negative
Beauville-Bogomolov square. Clearly, at least one of them necessarily has.

So the walls of K\"ahler chambers are parts of the orthogonal hyperplanes to the MBM classes, whereas, by Markman, the walls of the (larger) birational K\"ahler chambers are parts of the orthogonal hyperplanes to the classes of the prime uniruled divisors with negative square (``prime exceptional divisors'').
As we have seen in the previous section (proposition \ref{rat-quot}), every prime exceptional (in fact, even just uniruled) divisor is birationally a $\P^1$-bundle, and we may consider 
the class of its generating $\P^1$ in $H^2(X,\Q)$.

\begin{lemma} This class is MBM, 
proportional to the class of the divisor itself.
\end{lemma}

Indeed, both survive as $(1,1)$-classes on the same deformations of $X$.

\hfill

In this way, our concept of an MBM class is an extension of that of the prime exceptional divisor: in a suitable complex structure
$X=(M,I)$,
an MBM class is represented by a rational curve whose deformations within $X$ may cover a divisor or a lower-dimensional 
submanifold. This dimension is invariant by small deformations of $X$, see \cite{AV-MBM}. G. Mongardi in \cite{Mon} has introduced the notion of a {\bf wall divisor} which has turned out to be the same as our MBM class (but this was not immediately clear to us because of a somewhat misleading terminology, and also because the earlier versions of his paper dealt specifically with the K3 type case).

The birational K\"ahler cone $\Bir(X)$ is a connected component of the complement in $\Pos(X)$ to the union of hyperplanes 
$D^{\bot}$, where $D$ runs through the prime exceptional divisor classes. In other words the hyperplanes $D^{\bot}$ within $\Pos(X)$
are entirely made of the walls of birational K\"ahler chambers. Subdividing $\Bir(X)$ into K\"ahler cones of birational models of $X$ is more delicate: indeed a given rational curve $C$ may disappear on a given birational model, so it is not obvious
that the hyperplane $[C]^{\bot}$ continues through
$\Bir(X)$. We have however proved that this is true.

\begin{theorem} (\cite{AV-MBM}, Theorem 6.2)
The K\"ahler cone is a connected component of the
complement to the union of hyperplanes $z^\bot$ in
$\Pos(X)$, where $z$ runs through the set of MBM classes
of type $(1,1)$.
\end{theorem}

\begin{rem} It follows easily from the definition that the
  monodromy group preserves the set of MBM classes and
  that the Hodge monodromy  
	group\footnote{The Hodge monodromy
group $\Mon^{\Hdg}$ of an IHS manifold 
is the subgroup of the monodromy group 
(Subsection \ref{_Param_Subsection_}) 
preserving the Hodge decomposition.} $\Mon^{\Hdg}$
preserves the set of MBM classes of type $(1,1)$.
\end{rem}

\section{Results on MBM classes and applications}

\subsection{Markman's Torelli theorem and the birational
  cone conjecture}
\label{_Markman_Subsection_}

The following corollary of
 the global Torelli theorem 
due to Markman (see \cite{Mar-survey}, Th. 1.3, Cor. 5.7) is crucial for understanding the (birational) automorphisms
of an IHS manifold.

\begin{theorem} a) Let $\gamma\in O(H^2(X,\Z))$ be an element of the monodromy group preserving the Hodge decomposition and the K\"ahler cone.
	Then there exists $f\in \Aut(X)$ such that $f^*=\gamma$ (i.e. inducing $\gamma$ on cohomology).

	b) Let $\gamma\in O(H^2(X,\Z))$ be an element of the monodromy group preserving the Hodge decomposition and the birational K\"ahler cone.
	Then there exists a birational automorphism $f$ of $X$ such that $f^*=\gamma$.

\end{theorem}

In particular, in this way questions about automorphisms can be reduced to questions about the monodromy. With a little bit of linear algebra work (\cite{AV-MBM}, sections 3, 6) one sees that

\begin{itemize}
	\item To prove that $\Aut(X)$ acts with finitely
          many orbits on the walls of the K\"ahler cone,
          it suffices to prove that the Hodge monodromy group
$\Mon^{\Hdg}$ acts with finitely many orbits on the set of
          MBM classes of type $(1,1)$;

	\item To prove that the group of birational automorphisms acts with finitely many orbits on the walls of the birational K\"ahler cone, it suffices to prove that $\Mon^{\Hdg}$
		acts with finitely many orbits on the set of divisorial MBM classes of type $(1,1)$.

\end{itemize}

The second part has been settled by Markman quite some time before we have settled the first one in \cite{AV-cone} and was in fact a huge source of inspiration for all our work. Our point of view hopefully brought some simplifications to Markman's original arguments, however the crucial 
step (reflexions with respect to the divisorial MBM classes are in the monodromy) remains unchanged. The main steps of proof are as follows.

\medskip

{\bf Step 1}  We know that the Hodge monodromy is of finite index in $O(NS(X), q)$, therefore it suffices to prove that
$O(NS(X))$ acts with finitely many orbits on the set of prime exceptional divisor classes (divisorial MBM classes).

\smallskip

{\bf Step 2} By classical results on quadratic forms (\cite{Kn}, Satz 30.2, and a variant \cite{AV-MBM}, Theorem 3.14) this follows once we show that the divisorial
MBM classes have bounded square.

\smallskip

{\bf Step 3} Markman shows that for a divisorial MBM class $E$, the reflexion $R_E(x)=x-2q(E,x)/q(E,E)$ is in the (Hodge) monodromy group.

\smallskip

{\bf Step 4} By definition, the monodromy group is a subgroup of $O(H^2(X,\Z), q)$, in particular $R_E$ are integral. Hence $q(E,E)$ must be bounded in terms of the discriminant of $q$.

\medskip

The statement about the finiteness of the orbits is often
referred to as {\bf the weak Morrison-Kawamata cone
  conjecture}, and the statement about the rational
polyhedral fundamental domain in conjecture
\ref{cone-conj} as {\bf the strong Morrison-Kawamata cone
  conjecture}. The latter has a chance to be valid only in
the algebraic setting, i.e. for the ample/movable cone
rather than the K\"ahler/birational K\"ahler one. Indeed
for a very general IHS $X$, the K\"ahler cone is round
(it is equal to $\Pos(X)$)
whereas the group of (birational) automorphisms is trivial. 

In section 6 of \cite{Mar-survey} Markman has also proved the following birational version of the strong 
Morrison-Kawamata cone 
conjecture.

\begin{theorem} The group of birational automorphisms of $X$ has a rational polyhedral fundamental domain on the cone $\Mov^+(X)$, where $\Mov(X)$ is 
	the rational part of $\overline\Bir(X)$ (that is, the intersection of $\overline\Bir(X)$ and $NS(X)\otimes \R$) and $\Mov^+(X)$ stands for keeping only the rationally defined part of the boundary.
\end{theorem}

A few years later Markman and Yoshioka \cite{MY} have essentially proved the equivalence of the weak and the strong version for IHS, namely that the
boundedness of squares of primitive MBM classes implies the strong Morrison-Kawamata cone conjecture. Here is the original formulation of their result.

\begin{theorem} (\cite{MY}, theorem 1.3) Let $X$ be an IHS manifold. Assume that the Beauville-Bogomolov square of any integral primitive and extremal class in the Mori cone of any IHS birational model $Y$ is bounded from below by a constant, which depends only on the birational equivalence class of $X$. Then the strong Morrison-Kawamata conjecture holds for $X$.
\end{theorem}

As Markman and Yoshioka also remarked, the cone conjecture implies
that the number of isomorphism classes of 
IHS birational models of an IHS manifold is finite
(cf. \cite[Theorem 1.9]{AV-cone}).


\subsection{The cone conjecture via ergodic theory}

As we have mentioned already, the weak version of the cone conjecture is equivalent to the fact that the set of primitive  MBM classes of type $(1,1)$,
(or their orthogonal hyperplanes in $H^{1,1}(X,\R)$) is finite modulo the Hodge monodromy action. This in turn is equivalent
to the boundedness of the Beauville-Bogomolov square of
primitive MBM classes of type $(1,1)$ (Subsection 
\ref{_Markman_Subsection_}, Step 2). 
To prove this boundedness assertion, it suffices to do the case when 
$X$ has maximal Picard rank (indeed it is easy to see that we can deform $X$ to a manifold of maximal Picard rank in such a way that all 
MBM classes of type $(1,1)$ remain of type $(1,1)$). Such an $X$ satisfies $H^{1,1}(X,\R)=NS(X)\otimes \R$, the Beauville-Bogomolov form is 
of signature $(1,b_2-3)$ on $NS(X)$, in particular $X$ is projective (see \cite{Hu-basic}).  

It turns out that using ergodic theory one can prove the following general statement.

\begin{theorem} Let $V=V_\Z\otimes \R$ be a vector space with integral structure, equipped with an integral quadratic form $q$ of signature $(1,n)$. We suppose $n\geq 3$.
	Consider the projectivization $\P V^+$ of the positive cone of $q$, non-canonically isomorphic to $SO^+(1,n)/SO(n)$, as a model for the hyperbolic $n$-space $\H^n$. Let $\Gamma$ be a discrete subgroup of finite index in $SO^+(V_{\Z})$. Let $\Sigma$ be a $\Gamma$-invariant union of hyperplanes $S_i,\ i\in I$. Then either $\Sigma$ is dense in $\H^n$ or the set $I$
	is finite modulo  $\Gamma$.
\end{theorem}

The key point is that one can view the hyperplanes $S_i$ as the projections of orbits of a single subgroup 
$H \subset G=SO^+(1,n)$, $H\cong SO^+(1,n-1)$ (it is especially visible when $n=2$, though this case is excluded from the theorem: $G$ interprets as the unit tangent bundle to the hyperbolic plane $\H$, and by lifting the hyperplanes to $G$ tautologically we obtain a foliation, known as ``the geodesic flow''). When $n>2$, $H$ is generated by unipotents, which is the case
studied in Ratner's theory (Ratner's original papers are \cite{Ra1}, \cite{Ra2}; see also \cite{Morris} for an accessible exposition). Ratner gives a classification of $H$-invariant, $H$-ergodic measures on $G/\Gamma$, where $G$ is a non-compact simple Lie group with finite center, $\Gamma$ is a lattice (i.e. a discrete subgroup of finite index in $G$) and $H$ is a subgroup generated by unipotents. Her main result is that such a measure is {\bf algebraic}, i.e. is the pushforward of the
Haar measure of an orbit $Sx$ of a closed subgroup $S$, such that $x\Gamma x^{-1}$ is a lattice in $S$.

Moreover, a combination of results by Mozes-Shah and Dani-Margulis states that the weak limit of such measures either goes to infinity 
or is again a measure of the same type. The following formulation of these is Corollary 4.6 from \cite{AV-cone}.

\begin{theorem} Let $G$ be a connected Lie group, $\Gamma$ a lattice, ${\cal P}(Y)$ the space of
	all probability measures on $Y=G/\Gamma$, and ${\cal Q}(Y)\subset {\cal P}(Y)$  the space of all algebraic
probability measures. Then ${\cal Q}(Y)$  is closed in ${\cal P}(Y)$ with respect to weak topology. Moreover,
let $Y\cup\infty$ denote the one-point compactification of $Y$, so that ${\cal P}(Y\cup \infty)$ is compact. If
for a sequence $\mu_i\in {\cal Q}(Y)$, $\mu_i\to \mu \in {\cal P}(Y\cup \infty)$, then either $\mu\in {\cal Q}(Y)$,  or $\mu$ is supported
at infinity.
\end{theorem}

We associate measures $\mu_i$ to our hyperplanes $S_i$, by taking the closures of the corresponding orbits in $G/\Gamma$ (in fact they are closed already since the hyperplanes are rational). A purely geometric argument yields that no subsequence of  $\mu_i$ can tend to infinity (one easily reduces the problem to the case of geodesics on a hyperbolic Riemann surface, where it amounts to the statement that there are no closed geodesics around cusps). Any sequence of $\mu_i$ should then have a subsequence weakly converging to an algebraic measure. We remark that there is no intermediate closed subgroup between $H=SO^+(1,n-1)$ and $G=SO^+(1,n)$.  The case when the limit is supported on an orbit of $H$ amounts to finiteness (the converging subsequences must be constant in this case) and the case when the limit is supported on the whole of $G/\Gamma$ amounts to density. This is of course a very rough sketch; see \cite{AV-cone} for details.

\medskip

If $V=H^{1,1}(X,\R)$, $\Gamma$ is the Hodge monodromy and $\Sigma$ is the set of orthogonal hyperplanes to MBM classes, the first option of this alternative is clearly impossible, indeed the hyperplanes must bound the ample cone. Hence the cone conjecture holds:

\begin{theorem}\label{mbm-11}(\cite{AV-cone}) Let $X$ be an IHS manifold with $b_2(X)> 5$, then the Hodge monodromy group acts with finitely many orbits on the set of primitive MBM classes of type $(1,1)$. The Beauville-Bogomolov square of a primitive MBM class on $X$ is bounded by a constant $c=c(X)$.
\end{theorem}

\subsection{Uniform boundedness and an appication} 

The following theorem is a more precise version of Theorem
\ref{mbm-11} which also allows to settle the case $b_2=5$.

\begin{theorem}\label{mbm}(\cite{AV-orbits}) Let $X$ be an IHS manifold with $b_2(X)\geq 5$, then the monodromy group acts with finitely many orbits on the set of primitive MBM classes. The Beauville-Bogomolov square of a primitive MBM class on $X$ is bounded by a constant $c$ depending only on the deformation type of $X$ (or equivalently, on the underlying differentiable manifold $M$).
\end{theorem}

The idea behind the proof is similar, but one works with the homogeneous space $\Gr_{+++}(H^2(M,\R))$, isomorphic to $SO^+(3, b_2-3)/SO(3)\times SO(b_2-3)$, instead of the hyperbolic space $\H^{b_2-3}$ (as usually, $M$ denotes the underlying differentiable manifold).
There are two steps in the proof.

In the first step one observes that the subsets $z^{\bot}$ of $\Gr_{+++}(H^2(M,\R))$ 
(formed by positive 3-spaces orthogonal to an MBM class
$z$) are orbits of conjugated subgroups $P_z=\Stab (z)$ of
$G$ (we call them ``orbits of hyperplane type''). It turns
out that the orthogonality condition implies that all
these orbits are 
projections of right orbits of a single subgroup $P\cong
SO^+(3,b_2-4)$. $P$ is generated by unipotents when
$b_2\geq 5$, so Ratner theory and Mozes-Shah theorem apply
as soon as $b_2\geq 5$. One derives that the subsets
$z^{\bot}$ are either dense in $\Gr_{+++}(H^2(M,\R))$, or
finitely many up to the monodromy action. 

This part is independent of the hyperk\"aher context, just as Theorem 3.4. The geometric statement analogous to Theorem 3.4 which we prove is as follows. If $L$ is a lattice of signature $(p,q)$ with $p, q\geq 2$, $\Gamma$ is a finite index subgroup of $O(L)$, and $\Sigma$ an infinite union of $\Gamma$-orbits of negative vectors, then the set of positive $p$-planes orthogonal to a member of 
$\Sigma$ is dense in the positive grassmannian.

The second step is to show that the density is impossible in the hyperk\"ahler context, and therefore the set of the negative vectors
modulo $\Gamma$-action (that is, of the monodromy orbits of MBM classes) must be finite. To show this we remark that the complement 
to the union of all $z^{\bot}$ is identified with  
{\bf the Teichm\"uller space of hyperk\"ahler
  structures} on $M$, and therefore must have non-empty interior. This is the main content of the paper \cite{AV-teich}.
  In this paper we consider the space $\Hyp_m$ of all hyperk\"ahler metrics of volume 1 on a differentiable manifold $M$, and we
  define the Teichm\"uller space of hyperk\"ahler
  structures as the quotient of $\Hyp_m$ by isotopies.
  Each metric comes with a hyperk\"ahler triple $I,J,K$, inducing a triple of $2$-forms $\omega_I,\omega_J,\omega_K$.
  The corresponding triples of second cohomology classes generate positive three-planes, and it turns out that a positive three-plane corresponds to an actual hyperk\"ahler structure if it is not orthogonal to any MBM class; see  \cite{AV-teich} for details.

Note that Theorem \ref{mbm} is really a strengthening of
Theorem \ref{mbm-11}: if we know that the monodromy acts
with finitely many orbits on the set of primitive MBM
classes, we derive that the Beauville-Bogomolov square of
a primitive MBM class is bounded, in particular this holds
for MBM classes of type $(1,1)$ and so the Hodge monodromy
must act with finitely many orbits on the set of MBM
classes of type $(1,1)$.

\medskip

In this way it is easy to construct IHS without MBM classes of type $(1,1)$. It suffices to take a period point such that the Picard lattice of the corresponding complex structures does not represent small negative numbers (smaller than $c$ in absolute value). Of course one has to make sure that such points exist. But there are lattice-theoretic results on primitive embeddings from which one derives that any lattice of small rank $r$ and signature $(1,r-1)$ embeds primitively into 
$H^2(M, \Z)$.
Moreover, it follows from surjectivity of Torelli map that every primitive sublattice of $H^2(M, \Z)$ of signature $(1,r-1)$ is the Picard lattice of some complex structure (see \cite{Mor} for details in the case of K3 surfaces; the general IHSM case is completely analogous, \cite{AV-autom}).

An example of what one can prove in this way is the following result from \cite{AV-autom}.

\begin{theorem} Let $X$ be an IHS manifold with $b_2(X)\geq 5$. Then $X$ admits a projective deformation $X'$ with infinite group of symplectic automorphisms and Picard rank 2.
\end{theorem}

Indeed, it is possible to take a deformation $X'$ with Picard lattice $\Lambda$ of signature $(1,1)$, not representing
negative numbers of small absolute value nor zero. When $\Lambda$ does not represent zero, it is well-known that $O(\Lambda)$ is infinite. On the other hand, the Hodge monodromy group surjects onto a subgroup of finite index in $O(\Lambda)$, so it is infinite. Since there are no MBM classes of type $(1,1)$, the K\"ahler cone is equal to the positive cone and Markman's Hodge-theoretic Torelli theorem (Theorem 3.1 of the present survey) 
applies to conclude that every element of Hodge monodromy is induced by an automorphism.

\section{Contractibility and deformations}

For an MBM class $z$, the space
$\Teich_z$ is naturally split in two identical pieces, both
isomorphic to ${\cal P}_z=z^{\bot}\subset {\cal P}$ after
the Hausdorff reduction, according to whether $z$ 
takes positive or negative values on the K\"ahler cone.
This can be easily seen from Verbitsky's global Torelli theorem and the description 
of non-Hausdorff points in terms of
K\"ahler chambers (Markman's theorem 2.8).
 We call these two pieces
$\Teich^+_z$ resp. $\Teich^-_z$. Both are separated at a
general point of $\Teich_z$ (that is, where $z$ is the
only MBM class of type $(1,1)$) but not at points where
many other MBM classes are of type $(1,1)$. It is natural
to consider the subspace $\Teich_z^{min}$ of $\Teich^+_z$
where $z$ generates an extremal ray, in other words,
$z^{\bot}$ defines a wall of the K\"ahler cone: over each
period point where there are several K\"ahler chambers in 
$\Teich^+_z$, one keeps those which are adjacent to the
hyperplane $z^{\bot}$ and discards the rest. Along 
$\Teich^{min}_z$, a positive rational multiple of $z$ is
represented by an extremal rational curve.

By Kawamata-Shokurov base-point-free theorem, when $X$ is projective, extremal rational curves can be contracted. Indeed in this case the wall of the K\"ahler cone contained in $z^{\bot}$ has integral points. An integral point is the class of a line bundle $L$ which is nef and big. Since $X$ has trivial canonical class, Kawamata-Shokurov base-point-free theorem implies that some power $L^{\otimes n}$ is base-point-free. The morphism defined by the sections of $L^{\otimes n}$ is birational and contracts exactly the curves whose cohomology class is proportional to $z$.

A more recent result, essentially due to Bakker and Lehn, is that extremal rational curves also contract in the non-projective case.

\begin{theorem} Let $X$ be a hyperk\"ahler manifold with
  $b_2>5$ and $z$ an extremal class on $X$, that is, $z$
  is MBM and the complex structure $I$ of $X=(M,I)$ is in
  $\Teich^{min}_z$. Then there exists a bimeromorphic
  contraction $f:X\to Y$ contracting exactly the curves
   of cohomology class proportional to $z$.
\end{theorem}

The proof has two ingredients: local deformations and monodromy action. 

To describe the first one, take a projective IHS $X$ and
let $f:X\to Y$ be a birational morphism contracting the
class $z$. According to Namikawa (\cite{Nam} Section 3, p. 104) there is a
commutative diagram extending $f$:

$$\begin{CD}
{\cal X}  @>\Phi>>  {{\cal Y}}\\
@VVV              @VVV\\
 {\Def(X)}   @>G>>  {\Def(Y)}
\end{CD}$$

where $G$ is finite. This diagram in itself does not carry much information on contractions because the fiber 
of the family ${\cal Y}$ over a general $t\in \Def(Y)$ is smooth and the map $\Phi_t$ is an isomorphism even though 
$\Phi_0=f$. Bakker and Lehn have shown in \cite{BL} that $G$ induces an isomorphism between the subspace of ``locally trivial'' deformations $\Def^{lt}(Y)\subset \Def(Y)$ and the subspace $\Def_z(X)\subset \Def(X)$ of deformations preserving $z$ as a Hodge class. For $t\in \Def^{lt}(Y)$, $\Phi_t$ is a birational morphism contracting $z$. Here ``locally'' in the definition of local triviality means ``locally on $\cal Y$''
(every point of $\cal Y$ has a small neighbourhood which decomposes into a product, see \cite{FK} for precise definition and more), so this is an equisingularity condition.

The upshot is that if $z$ can be contracted on a projective $X$, then it also can be contracted on its sufficiently close non-projective neighbours. 

The second ingredient is as follows. On each connected
component of $\Teich^{min}_z$, there is an action of the
group $\Gamma_z$, defined as the subgroup
of the mapping class group preserving the component and
the class $z$. One shows exactly in the same way as in
\cite{V-erg}, \cite{V-erg-err},
that there are three types of $\Gamma_z$-orbits according to the rationality properties of the period plane. Here it is important to assume that $b_2(X)>5$ (Verbitsky's description is valid for $b_2(X)\geq 5$, but ${\cal P}_z$ and $\Teich_z$ are hyperplanes in ${\cal P}$ resp. in $\Teich$).
So the local picture globalizes as follows. For a non-projective $X'$ we can find a projective $X$ such that the monodromy transform of $X'$
(by some element of $\Gamma_z$) lands in an arbitrarily small neighbourhood of $X$: we can take any $X$ if the orbit of $X'$ is dense, and a projective $X$ from the closure of a non-dense orbit otherwise. So $z$ contracts on this monodromy transform, and consequently on $X'$ itself, since the action of the monodromy is the transport of complex structure which moreover preserves $z$.

\begin{rem} This result yields some information on the structure of the
rational quotient map already mentioned in 
Proposition \ref{rat-quot}. Namely, the rational quotient of any
component of the exceptional locus $Z$ is a regular map.
Indeed, this is the restriction of the contraction map to
$Z$.
\end{rem}

In \cite{AV-contr}, we prove that Bakker and Lehn's local triviality implies the local triviality in the usual sense of the word (local on the base, that is, any fiber $Y_t$ of the family ${\cal Y}^{lt}$ has a neighbourhood which decomposes into a product) in the  
{\bf real analytic category}. This implies in particular that the contraction loci are diffeomorphic as we move $X$ along $\Teich^{min}_z$.
Using ergodic actions, one can prove, with some restrictions, a stronger result. Namely, the fibers of the rational quotient map keep their biholomorphism type as $X$ varies in $\Teich^{min}_z$, under the condition that they are normal and assuming that the Picard number of $X$ is not maximal. See \cite{AV-contr} for details.

\begin{rem} As an obvious example of the exceptional divisor of the Hilbert square of a K3 surface $S$ shows, one cannot expect that the entire contraction loci keep their biholomorphism type. Indeed this exceptional divisor is a $\P^1$-bundle over $S$, so its biholomorphism type does vary with $S$.
\end{rem}

\section{Classification of MBM classes in low dimension for K3 type}

On $K3$ surfaces, primitive MBM classes are exactly the classes of square $-2$ in $H^2(M, \Z)$. On IHS fourfolds of $K3$ type, there are three ``obvious'' types of extremal curves/MBM classes, described by Hassett and Tschinkel in \cite{HT1}. Indeed first of all there are ``exceptional'' projective lines contracted by the Hilbert-Chow map. Since the restriction of the exceptional divisor to such a curve has degree $-2$ (blow-up of double points), one computes that the class of such a line $l_1$ is half-integral and its Beauville-Bogomolov square equals $-1/2$. To see the other two,
take a K3 surface $Y$ with a $(-2)$-curve $C$ and let $X$ be the Hilbert square of $Y$. Then $X$ carries at least two MBM classes: that of a line $l_2$ in the lagrangian plane $S^2C\subset S^{[2]}Y$ (where $S^2$ means the symmetric square) and that of $C$ itself viewed as a ruling of the image of $C\times Y$. The class $C$ is integral of square $-2$, and $l_2$ is half-integral of square $-5/2$ (see e.g. \cite{HT2}). Note that the class of $C$ is not extremal: contracting $C$, one automatically contracts $l_2$.
However it is MBM, and in fact becomes extremal when we perform the Mukai flop of $X$ centered at $S^2C$.

Using the deformation-invariance properties of MBM classes
and MBM loci, we provide an easy proof that there are no
other MBM classes (up to monodromy) in \cite{AV-K3}. The
key point of our argument is that, given an IHS fourfold
$X$ of K3 type and a negative integral $(1,1)$-class
$\eta$ on $X$, 
we can deform $X$, keeping $\eta$ of type $(1,1)$, to an
IHS $X'$ with negative definite Picard lattice of rank
two, which also happens to be the Hilbert scheme of a
nonprojective K3 surface. 
Indeed, one can detect
whether a given IHS manifold of K3 type is a 
Hilbert scheme of a K3 surface from its periods. It is a 
Hilbert scheme if it is a monodromy transform of a
manifold $X$ such that $e$ is of type (1,1),
where  $e$ is, as usual, one half
of the class of the exceptional divisor (under some fixed identification of $H^2(M,\Z)$ with the second cohomology 
of a Hilbert scheme).
Therefore, the periods of Hilbert schemes
are a union of countably many divisors
in $\Per$ realized as $\gamma(e)^{\bot}$, for
all $\gamma$ in the monodromy group. 

It turns out that one can find an
element $\gamma$ in the monodromy such that the sublattice
in the second cohomologies generated by $\eta$ and
$\gamma(e)$ is negative definite. Then for $X'$ as above, one can take a general complex
structure in the intersection of $\Teich_{\eta}$ and
$\Teich_{\gamma(e)}$.

On $X'$ there are not many rational curves: indeed $X'$ is the Hilbert square of a K3 surface $Y'$ which contains a single $(-2)$-curve and no other curves. All MBM classes and their contraction loci are therefore as described in the previous paragraph. Clearly $\eta$ must be one of them, and if $\eta$ is extremal, its contraction locus is either the projective plane or a $\P^1$-bundle over a K3 surface. The same is true on $X$ by deformation invariance properties. 

An interesting phenomenon occurs when we study the MBM classes on an IHS sixfold $X$ of K3 type. It turns out that one can not always 
deform $X$, preserving $\eta$, to the Hilbert cube of a K3 surface with negative cyclic Picard group. However, one can always deform it 
to the  Hilbert cube of a K3 surface with non-positive cyclic Picard group. On such manifolds, it is still easy to describe all rational 
curves. In this way, one gets a ``new'' MBM class (with respect to the four types obtained from deformation to the Hilbert cube of a 
K3 surface with Picard group generated by a $(-2)$-curve, as we did for Hilbert square). This is also the class which was missing from Hassett and Tschinkel's tables in \cite{HT2} and which is present in their later elaboration of Bayer and Macri's computations for low-dimensional IHS of K3 type  \cite{HT3}.

We refer the reader to \cite{AV-K3} for details and further analysis.

\section{Some open questions}

All known examples of IHS manifolds carry MBM classes. Can one prove that this must always be the case?

Markman's result that the reflections in divisorial MBM classes are elements of the monodromy group is extremely important.
One may ask how the non-divisorial MBM classes contribute to the monodromy. If a reflexion in a non-divisorial MBM classes is (or is not) in the monodromy, what geometric information does this give?

Though one has a good idea about the geometry of subvarieties covered by extremal rational curves, many details 
are still to be verified. For example
a Lagrangian subvariety covered by extremal rational curves tends to be isomorphic to the projective space. 
Wisniewski and Wierzba proved that on a holomorphic symplectic fourfold it is actually a plane (\cite{WW}). 
In higher dimension one can observe that for a Lagrangian submanifold $Y\subset X$, the normal bundle $N_{Y,X}$ is isomorphic to the cotangent $\Omega^1_Y$ and speculate that a subvariety contractible to a point 
should have ``negative'' normal bundle, so that $\Omega^1_Y$ is ``negative'' and $T_Y$ is ``positive''. If ``positive'' meant ample, then we would conclude that
$Y$ is $\P^n$, at least when smooth (\cite{Mo}). This however does not
work literally, because the contractibility of $Y$ does not mean that $N^*_{Y,X}$ is ample unless $Y$ is of codimension one. By a result of Ancona and Vo Van Tan, the contractibility of $Y$ implies that there exists a scheme structure $Y'$ on $Y$ such that the conormal sheaf of $Y'$ is ample (see \cite{AT} for details; we are grateful to Andreas H\"oring for indicating this result and the reference). In their setting,
$Y$ does not have to be smooth. The question about the possible isomorphism type of contractible Lagrangian subvarieties remains open.

{\scriptsize

\noindent {\sc Ekaterina Amerik\\
{\sc Universit\'e Paris-Saclay,\\
Laboratoire de Math\'ematiques d'Orsay,\\
Campus d'Orsay, B\^atiment 307, 91405 Orsay, France},\\
also:\\
{\sc Laboratory of Algebraic Geometry,\\
National Research University HSE,\\
Department of Mathematics, 6 Usacheva Str. Moscow, Russia,}\\
\tt  Ekaterina.Amerik@gmail.com},

\hfill

\noindent {\sc Misha Verbitsky\\
{\sc Instituto Nacional de Matem\'atica Pura e
              Aplicada (IMPA) \\ Estrada Dona Castorina, 110\\
Jardim Bot\^anico, CEP 22460-320\\
Rio de Janeiro, RJ - Brasil }\\
also:\\
{\sc Laboratory of Algebraic Geometry,\\
National Research University HSE,\\
Department of Mathematics, 6 Usacheva Str. Moscow, Russia,}\\
\tt  verbit@impa.br}.
 }


\begin{thebibliography}{} 
 
\bibitem[AV1]{AV-MBM}  E. Amerik, M. Verbitsky: Rational curves on hyperkähler manifolds, Int.
Math. Res. Notices 2015, no.23, 13009 -- 13045

\bibitem[AV2]{AV-cone} 
E. Amerik, M. Verbitsky: Morrison-Kawamata cone conjecture
for hyperk\"ahler manifolds. Ann. Sci. ENS (4) 50 (2017),
no. 4, 973-993


\bibitem[AV3]{AV-orbits} E. Amerik, M. Verbitsky: Collections of Orbits of Hyperplane Type in Homogeneous Spaces, Homogeneous Dynamics, and Hyperk\"ahler Geometry, IMRN 2020, no. 1 (2020), 25--38.

\bibitem[AV4]{AV-teich} E. Amerik, M. Verbitsky: Teichm\"uller space for hyperk\"ahler and symplectic structures. J. Geom. Phys. 97 (2015), 44--50.

\bibitem[AV5]{AV-autom} E. Amerik, M. Verbitsky: Construction of automorphisms of hyperk\"ahler manifolds. Compos. Math. 153 (2017), no. 8, 1610--1621.

\bibitem[AV6]{AV-contr} E. Amerik, M. Verbitsky: Contraction centers in families of hyperk\"ahler manifolds, https://arxiv.org/abs/1903.04884

\bibitem[AV7]{AV-K3} E. Amerik, M. Verbitsky: MBM classes and contraction loci on low-dimensional hyperk\"ahler manifolds of K3 type, 
https://arxiv.org/abs/1907.13256

\bibitem[AT]{AT} V. Ancona, Vo Van Tan: On the blowing down problem in C-analytic geometry. J. Reine Angew. Math. 350 (1984), 178--182

\bibitem[Be]{Be} A. Beauville, Vari\'et\'es K\"ahleriennes dont la premi\`ere classe de Chern est nulle.
J. Differential Geom. 18 (1983), no. 4, 755–782 (1984)

\bibitem[BL]{BL} B. Bakker, C. Lehn: A global Torelli theorem for singular symplectic varieties, https://arxiv.org/abs/1612.07894

\bibitem[BM]{BM} A. Bayer, E. Macri: MMP for moduli of sheaves on K3s via wall-crossing: nef and movable cones, Lagrangian fibrations. Invent. Math. 198 (2014), no. 3, 505--590. 

\bibitem[BHT]{BHT} A. Bayer, B. Hassett, Y. Tschinkel: Mori cones of holomorphic symplectic varieties of K3 type. Ann. Sci. Éc. Norm. Supér. (4) 48 (2015), no. 4, 941--950.

\bibitem[Bou]{Bou} 
S. Boucksom, Divisorial Zariski decompositions on compact complex manifolds,
Ann. Sci. Ecole Norm. Sup. (4) 37 (2004), no. 1, 45--76

\bibitem[D]{D} S. Druel, Quelques remarques sur la d\'ecomposition de Zariski divisorielle sur les vari\'et\'es dont la premi\`ere classe de Chern est nulle. Math. Z. 267 (2011), no. 1-2, 413--423. 

\bibitem[FK]{FK} H. Flenner, S. Kosarew, On locally trivial deformations.
Publ. Res. Inst. Math. Sci. 23 (1987), no. 4, 627--665. 

\bibitem[F] {F} A. Fujiki, On the de Rham cohomology group of a compact K\"ahler symplectic manifold, in: Advanced Studies in Pure Mathematics 10 (1987)

\bibitem[HT1]{HT1} B. Hassett, Yu. Tschinkel, Rational curves on symplectic fourfold, Geom. Funct. Anal. 11 (2001), no. 6, 1201--1228.

\bibitem[HT2]{HT2} B. Hassett, Yu. Tschinkel, Intersection numbers of extremal rays on holomorphic
symplectic varieties. Asian J. Math. 14 (2010), no. 3, 303--322

\bibitem[HT3]{HT3} B. Hassett, Yu. Tschinkel, Extremal rays and automorphisms of holomorphic
symplectic varieties. K3 surfaces and their moduli, 73--95, Progr. Math.,
315, Birkh\"auser/Springer, [Cham], 2016.


\bibitem[Hu1]{Hu-basic} D. Huybrechts, Compact
  hyper-K\"ahler manifolds: basic results, erratum,
  Invent. Math. 152 (2003), no. 1, 209-212. 

\bibitem[Hu2]{_Huybrechts:cone_}
Huybrechts, D.,
{\em The K\"ahler cone of a compact hyperk\"ahler manifold},
Math. Ann. 326 (2003), no. 3, 499--513, arXiv:math/9909109.

\bibitem[Ka]{_Kawamata:cone_} 
Y. Kawamata,
{\em On the cone of divisors of Calabi-Yau fiber spaces}, 
Int. J. Math.
     8 (1997), 665--687.

\bibitem[KM]{KM} J. Kollar, S. Mori, Birational geometry of algebraic varieties, Cambridge Tracts in Mathematics 134, Cambridge Univ. Press. 


\bibitem[Kn]{Kn} M. Kneser, R. Scharlau, Quadratische formen, Berlin, Springer, 2002.

\bibitem[KS]{_Kod-Spen-AnnMath-1960_}
Kodaira, K., Spencer, D.C., {\em On deformations of complex analytic structures. III. Stab
ility theorems for complex structures},
Ann. of Math. \textbf{71} (1960), 43-76.


\bibitem[LP]{LP} E. Looienga, C. Peters Torelli theorems for K\"ahler K3 surfaces.
Compositio Math. 42 (1980/81), no. 2, 145--186. 


\bibitem[Mar]{Mar-survey} E. Markman, A survey of Torelli and monodromy results for holomorphic symplectic varieties, in Complex and differential geometry, Springer Proc. Math. 8, Springer, Heidelberg, 2011, 257--322 

\bibitem[MY] E. Markman, K. Yoshioka, A proof of the Kawamata-Morrison cone conjecture for holomorphic symplectic varieties of K3 or generalized Kummer deformation type,
Int. Math. Res. Not. IMRN 2015, no. 24, 13563--13574. 

\bibitem[Mon]{Mon} G. Mongardi, A note on the K\"ahler and Mori cones of hyperk\"ahler manifolds, Asian J. Math. 19 (2015), 583--591. 

\bibitem[Mo]{Mo} S. Mori, Projective manifolds with ample tangent bundles. Ann. of Math. (2) 110 (1979), no. 3, 593--606

\bibitem[Mor1]{Mor} D. Morrison, On K3 surfaces with large Picard number. Invent. Math. 75 (1984), no. 1, 105--121.

\bibitem[Mor2]{_Morrison:Beyond_}
D. Morrison,
{\em Beyond the K\"ahler cone}, 
Proceedings of the Hirzebruch 65 
conference on algebraic geometry (Ramat Gan, 1993), 
361-376, Bar-Ilan Univ.  (1996).



\bibitem[Morr]{Morris} D. W. Morris, Ratner's theorems on unipotent flows, Chicago Lectures in Mathematics, University of Chicago Press, Chicago, IL, 2005.

\bibitem[MY]{MY} E. Markman, K. Yoshioka: A proof of the Kawamata-Morrison cone conjecture for holomorphic symplectic varieties of K3${}^{[n]}$ or generalized Kummer deformation type,
Int. Math. Res. Not. IMRN 2015, no. 24, 13563--13574. 

\bibitem[N]{Nam} Y. Namikawa, On deformations of Q-factorial symplectic varieties. J. Reine Angew. Math. 599 (2006), 97--110.

\bibitem[PS-S]{PS-S} I. Pyatetski-Shapiro, I. Shafarevich: Torelli's theorem for algebraic surfaces of type K3. (Russian) Izv. Akad. Nauk SSSR Ser. Mat. 35 1971 530--572

\bibitem[Ran]{Ran} Z. Ran, Hodge theory and deformations of maps, Compositio Math. 97 (1995), no. 3, 309--328.

\bibitem[Ra1]{Ra1} M. Ratner, On Raghunathan's measure conjecture, Ann. of Math. 134 (1991), 545--607

\bibitem[Ra2]{Ra2} M. Ratner, Raghunathan's topological conjecture and distributions of unipotent flows, Duke Math. J. 63 (1991), 235--280. 

\bibitem[St]{St} H. Sterk, Finiteness results for algebraic K3 surfaces. Math. Z. 189 (1985), no. 4, 507--513.

\bibitem[T]{_Totaro:MK_cone_}
Burt Totaro, 
{\em The cone conjecture for Calabi-Yau pairs in dimension
  2},  Duke Math. J. Volume 154, Number 2 (2010), 241-263. 


\bibitem[V1]{V-Tor} M. Verbitsky, Mapping class group and a global Torelli theorem for hyperk\"ahler manifolds. Appendix A by Eyal Markman. 
	Duke Math. J. 162 (2013), no. 15, 2929--2986. 

\bibitem[V1bis]{V-Tor-err}  M. Verbitsky: Mapping class group and a global Torelli theorem, erratum, arxiv https://arxiv.org/abs/1908.11772

\bibitem[V2]{V-erg} 
M. Verbitsky, Ergodic complex structures on hyperk\"ahler manifolds. Acta Math. 215 (2015), no. 1, 161--182.

\bibitem[V2bis]{V-erg-err} M. Verbitsky, Ergodic complex structures on hyperk\"ahler manifolds, erratum, https://arxiv.org/abs/1708.05802

\bibitem[WW]{WW} Wierzba, Jan; Wisniewski, Jaroslaw A. Small contractions of symplectic 4-folds. Duke Math. J. 120 (2003), no. 1, 65--95.

\end{thebibliography}
\end{document}